\documentclass[11pt]{article}

\usepackage{amsmath,amssymb}
\usepackage{latexsym}
\usepackage{graphicx}
\usepackage{bm}

\newtheorem{thm}{Theorem}

\newtheorem{lem}[thm]{Lemma}
\newtheorem{rem}[thm]{Remark}

\newenvironment{proof}{\begin{trivlist}
                       \item[]{\bf Proof.}
                       \hspace{0cm}}{\hfill $\Box$
                       \end{trivlist}}

\begin{document}
\title{An iterative scheme for solving nonlinear equations\\ with monotone operators}

\author{N. S. Hoang$\dag$\footnotemark[1]
 \quad 	  A. G. Ramm$\dag$\footnotemark[3] \\
\\
\\
$\dag$Mathematics Department, Kansas State University,\\
Manhattan, KS 66506-2602, USA
}

\renewcommand{\thefootnote}{\fnsymbol{footnote}}
\footnotetext[1]{Email: nguyenhs@math.ksu.edu}
\footnotetext[3]{Corresponding author. Email: ramm@math.ksu.edu}

\date{}
\maketitle

\begin{abstract}
\noindent
An iterative scheme for solving ill-posed
nonlinear operator equations with monotone operators
is introduced and studied in this paper. A Dynamical Systems Method
(DSM) algorithm for stable solution of ill-posed operator
equations with monotone operators is proposed and its convergence is proved.
A new discrepancy principle is proposed and justified.
A priori and a posteriori stopping rules for the DSM algorithm
are formulated and justified.
\\


\end{abstract}

\section{Introduction}

In this paper we study a Dynamical Systems Method (DSM) for solving the equation
\begin{equation}
\label{aeq1}
F(u)=f,
\end{equation}
 where
$F$ is a nonlinear twice Fr\'{e}chet differentiable monotone operator in a 
real Hilbert space $H$, and equation \eqref{aeq1} is assumed solvable. 
Monotonicity is understood in the following sense: 
\begin{equation}
\label{ceq2}
\langle F(u)-F(v),u-v\rangle\ge 0,\quad \forall u,v\in H.
\end{equation}
It is known (see, e.g., \cite{R499}), that the set $\mathcal{N}:=\{u:F(u)=f\}$ is closed and convex if $F$ is monotone and continuous. 
A closed and convex set in a Hilbert space has a unique minimal-norm element. 
This element in $\mathcal{N}$ we denote by $y$, $F(y)=f$. 
We assume that 
\begin{equation}
\label{ceq3}
\sup_{\|u-u_0\|\le R}\|F^{(j)}(u)\|\le M_j(R),\quad 0\le j\le 2,
\end{equation}
where $u_0\in H$ is an element of $H$,
$R>0$ is arbitrary, 
and $f=F(y)$ is not known but
 $f_\delta$, the noisy data, are known and $\|f_\delta-f\|\le \delta$. If $F'(u)$
 is not boundedly invertible then 
 solving for $u$ given noisy data $f_\delta$ is often (but not always) an ill-posed problem.

Our goal is to develop an iterative process for a stable solution of equation \eqref{aeq1},
given noisy data $f_\delta$, $\|f-f_\delta\|\le \delta$. 
The idea of this iterative process is similar to the ideas of the DSM method in \cite{R540}, 
\cite{R499}--\cite{R491}. In \cite{I}--\cite{R499} and references therein methods for solving 
ill-posed problems are discussed.

\section{Auxiliary and main results}

We assume throughout the paper that $0<(a_n)_{n=0}^\infty \searrow 0$.

\subsection{Auxiliary results}

The inner product in $H$ is denoted $\langle u,v\rangle$. 
Let us consider the following equation
\begin{equation}
\label{21eq2}
F(V_\delta)+a V_\delta-f_\delta = 0,\qquad a>0.
\end{equation}
It is known (see, e.g., \cite{D}) that equation \eqref{2eq2} with monotone continuous operator $F$ has a unique solution for any 
fixed $a>0$ and $f_\delta\in H$. 

\begin{lem}
\label{lem0}
If \eqref{ceq2} holds and $F$ is continuous, then
$\|V_\delta\|=O(\frac{1}{a})$ as $a\to\infty$, and
\begin{equation}
\label{4eq2}
\lim_{a\to\infty}\|F(V_\delta)-f_\delta\|=\|F(0)-f_\delta\|.
\end{equation}
\end{lem}

\begin{proof}
Rewrite \eqref{2eq2} as
$$
F(V_\delta) - F(0) + aV_\delta + F(0)-f_\delta = 0.
$$
Multiply this equation
by $V_\delta$, use the inequality $\langle F(V_\delta)-F(0),V_\delta-0\rangle \ge 0$, which
follows from \eqref{ceq2}, and get:
$$
a\|V_\delta\|^2\le\langle aV_\delta + F(V_\delta)-F(0), V_\delta\rangle = 
\langle f_\delta-F(0), V_\delta\rangle \le \|f_\delta-F(0)\|\|V_\delta\|.
$$
Therefore, 
$\|V_\delta\|=O(\frac{1}{a})$. This and the continuity of $F$ imply \eqref{4eq2}.
\end{proof}

Let us consider the following equation
\begin{equation}
\label{2eq2}
F(V_{n,\delta})+a_n V_{n,\delta}-f_\delta = 0,\qquad a_n>0.
\end{equation}
Let us denote $V_n:=V_{n,\delta}$ when $\delta\not=0$. 
From the triangle inequality one gets:
$$
\|F(V_0)-f_\delta\|\ge\|F(0) - f_\delta\| -\|F(V_0)-F(0)\|.
$$
From the inequality $\|F(V_0)-F(0)\|\le M_1\|V_0\|$ and Lemma~\ref{lem0} it follows that for large $a_0$ one has:
$$
\|F(V_0)-F(0)\|\le M_1\|V_0\|=O\bigg{(}\frac{1}{a_0}\bigg{)}.
$$
Therefore, 
if $\|F(0)-f_\delta\|>C\delta$, then $\|F(V_0)-f_\delta\|\ge (C-\epsilon)\delta$, 
where $\epsilon>0$ is arbitrarily small, for sufficiently large $a_0>0$. 

\begin{lem}
\label{lem11}
Suppose that $\|F(0)-f_\delta\|> C\delta$, \,$C>1$. 
Assume that $0<(a_n)_{n=0}^\infty \searrow 0$, and $a_0$ is sufficiently large. 
Then, there exists a unique $n_\delta>0$, such that 
\begin{equation}
\label{eeq5}
\|F(V_{n_\delta})-f_\delta\|\le C\delta < \|F(V_n)-f_\delta\|, \quad \forall n<n_\delta.
\end{equation}
\end{lem}

\begin{proof}
We have $F(y)=f$, and
\begin{align*}
0&=\langle F(V_n)+a_nV_n-f_\delta, F(V_n)-f_\delta \rangle\\
&=\|F(V_n)-f_\delta\|^2+a_n\langle V_n-y, F(V_n)-f_\delta \rangle + a_n\langle y, F(V_n)-f_\delta \rangle\\
&=\|F(V_n)-f_\delta\|^2+a_n\langle V_n-y, F(V_n)-F(y) \rangle + a_n\langle V_n-y, f-f_\delta \rangle 
+ a_n\langle y, F(V_n)-f_\delta \rangle\\
&\ge\|F(V_n)-f_\delta\|^2 + a_n\langle V_n-y, f-f_\delta \rangle + a_n\langle y, F(V_n)-f_\delta \rangle.
\end{align*}
Here the inequality $\langle V_n-y, F(V_n)-F(y) \rangle\ge0$ was used. 
Therefore
\begin{equation}
\label{1eq1}
\begin{split}
\|F(V_n)-f_\delta\|^2 &\le -a_n\langle V_n-y, f-f_\delta \rangle - a_n\langle y, F(V_n)-f_\delta \rangle\\
&\le a_n\|V_n-y\| \|f-f_\delta\| + a_n\|y\| \|F(V_n)-f_\delta\|\\
&\le  a_n\delta \|V_n-y\|  + a_n\|y\| \|F(V_n)-f_\delta\|.
\end{split}
\end{equation}
On the other hand, one has:
\begin{align*}
0&= \langle F(V_n)-F(y) + a_nV_n +f -f_\delta, V_n-y\rangle\\
&=\langle F(V_n)-F(y),V_n-y\rangle + a_n\| V_n-y\| ^2 + a_n\langle y, V_n-y\rangle + \langle f-f_\delta, V_n-y\rangle\\
&\ge  a_n\| V_n-y\| ^2 + a_n\langle y, V_n-y\rangle + \langle f-f_\delta, V_n-y\rangle,
\end{align*}
where the inequality $\langle V_n-y, F(V_n)-F(y) \rangle\ge0$ was used. Therefore,
$$
a_n\|V_n-y\|^2 \le a_n\|y\|\|V_n-y\|+\delta\|V_n-y\|.
$$
This implies
\begin{equation}
\label{1eq2}
a_n\|V_n-y\|\le a_n\|y\|+\delta.
\end{equation}
From \eqref{1eq1} and \eqref{1eq2}, and an elementary inequality $ab\le \epsilon a_n^2+\frac{b^2}{4\epsilon},\,\forall\epsilon>0$, one gets:
\begin{equation}
\label{3eq4}
\begin{split}
\|F(V_n)-f_\delta\|^2&\le \delta^2 + a_n\|y\|\delta + a_n\|y\| \|F(V_n)-f_\delta\|\\
&\le \delta^2 + a_n\|y\|\delta + \epsilon \|F(V_n)-f_\delta\|^2 + 
\frac{1}{4\epsilon}a_n^2\|y\|^2,
\end{split}
\end{equation}
where $\epsilon>0$ is fixed, independent of $n$, and can be chosen arbitrary small. 
Let $n\to\infty$ so $a_n\searrow 0$. Then \eqref{3eq4} implies
$\lim_{n\to\infty}(1-\epsilon)\|F(V_n)-f_\delta\|^2\le \delta^2$,\, $\forall\, \epsilon>0$. This implies
$\lim_{n\to\infty}\|F(V_n)-f_\delta\| \le \delta$.
This, the assumption $\|F(0)-f_\delta\|>C\delta$, 
and the fact that $\|F(V_n)-f_\delta\|$ is nonincreasing (see Lemma~\ref{remark1}),
imply that there exists a unique $n_\delta>0$ such that \eqref{eeq5} holds. 
Lemma~\ref{lem11} is proved.
\end{proof}

\begin{rem}
\label{rem3}
{\rm 
Let $V_{0,n}:=V_{\delta,n}|_{\delta=0}$. Then $F(V_{0,n})+a_n V_{0,n}-f=0$. 
It was proved in \cite{R499} that $\|V_{0,n}\|\le \|y\|$, where $y$
is the minimal-norm solution to $F(u)=f$, and that 
$$
\|V_{\delta,n}-V_{0,n}\|\le \frac{\delta}{a_n}.
$$
Thus, one gets the following estimate:
\begin{equation}
\label{2eq1}
\|V_{n}\|\le \|V_{0,n}\|+\frac{\delta}{a_n}\le \|y\|+\frac{\delta}{a_n}, 
\quad V_n:=V_{\delta,n}.
\end{equation}
}
\end{rem}

\begin{lem}
\label{remark1}
Let $0<a_n\searrow 0$, and $F$ be monotone.
Denote 
$$
h_n:=\|F(V_n) - f_\delta\|,\quad  g_n:=\|V_n\|,\qquad n=0,1,...,
$$ 
where $V_n$ solves \eqref{2eq2} with $a=a_n$. 
Then
$h_n$ is nonincreasing, and $g_n$ is nondecreasing.
\end{lem}
 
\begin{proof}
Note that $h_n=a_n\|V_n\|$. One has
\begin{equation}
\label{1eq3}
\begin{split}
0&\le \langle F(V_n)-F(V_m),V_n-V_m\rangle\\
&= \langle -a_nV_n + a_mV_m,V_n-V_m\rangle\\
&= (a_n+ a_m)\langle V_n,V_m \rangle -a_n\|V_n\|^2 -  a_m\|V_m\|^2.
\end{split}
\end{equation}
Thus,
\begin{equation}
\label{2eq6}
\begin{split}
0&\le (a_n+ a_m)\langle V_n,V_m \rangle -a_n\|V_n\|^2 -  a_m\|V_m\|^2\\
& \le  (a_n+ a_m)\|V_n\|\|V_m \| - a_n\|V_n\|^2 -  a_m\|V_m\|^2\\
& = (a_n \|V_n\| -  a_m \|V_m\|)(\|V_m\|-\|V_n\|)\\
& = (h_n-h_m)(g_m-g_n).
\end{split}
\end{equation}
If $g_m\ge g_n$ then \eqref{2eq6} implies $h_n\ge h_m$, so
$$
a_ng_n\ge  a_mg_m\ge  a_m g_n.
$$
Thus, if $g_m\ge g_n$ then $ a_m\le a_n$ and, therefore, $m\ge n$,
because $a_n$ is decreasing.

Conversely, if $g_m\le g_n$ then $h_n\le h_m$.
This implies $ a_m\ge a_n$, so $m\le n$.

Therefore $h_n$ is nonincreasing
and $g_n$ is nondecreasing. Lemma~\ref{remark1} is proved.
\end{proof}

\begin{rem}
\label{remmoi}
{\rm
From Lemma~\ref{lem0} and Lemma~\ref{remark1} one concludes that
$$
a_n\|V_n\|=\|F(V_n)-f_\delta\| \le \|F(0)-f_\delta\|,\qquad \forall n\ge 0.
$$
}
\end{rem}

\begin{lem}
\label{lem1}
Suppose $M_1, c_0$, and $c_1$ are positive constants and $0\not=y\in H$.
Then there exist $\lambda>0$ and a sequence $0<(a_n)_{n=0}^\infty\searrow 0$ such that the following conditions hold
\begin{align}
\label{yeq22}
a_n &\le 2a_{n+1},\\
\label{yeq23}
\|f_\delta -F(0)\| &\le \frac{a_0^2}{\lambda},\\
\label{yeq24}
\frac{M_1}{\lambda}  &\le \|y\|,\\
\label{yeq25}
\frac{a_n - a_{n+1}}{a_{n+1}^2} &\le \frac{1}{2c_1\lambda},\\
\label{yeq26}
c_0\frac{a_n}{\lambda^2}+ \frac{a_n-a_{n+1}}{a_{n+1}}c_1 &\le \frac{a_{n+1}}{\lambda}.
\end{align}
\end{lem}

\begin{proof}
Let us show that if $0<a_0$ is sufficiently large then the following sequence
\begin{equation}
\label{peq19}
a_n = \frac{a_0}{1+n},
\end{equation}
satisfy conditions \eqref{yeq22}--\eqref{yeq26}. 
One has
$$
\frac{a_n}{a_{n+1}} = \frac{n+2}{n+1} \le 2,\qquad \forall\, n\ge 0.
$$
Thus, inequality \eqref{yeq22} is obtained.

Choose 
\begin{equation}
\label{peq20}
\lambda \ge \frac{M_1}{\|y\|}
\end{equation}
then inequality \eqref{yeq24} is satisfied. 

Inequality \eqref{yeq23} is obtained if $a_0$ is sufficiently large. Indeed, \eqref{yeq23} holds
if
\begin{equation}
\label{qeq21}
a_0 \ge \sqrt{\lambda\|f_\delta -F(0)\|}.
\end{equation}

Let us check inequality \eqref{yeq25}. One has
$$
\frac{a_n-a_{n+1}}{a_{n+1}^2} = \bigg{(}\frac{a_0}{1+n}-\frac{a_0}{2+n}\bigg{)}\frac{(n+2)^2}{a_0^2}
=\frac{n+2}{a_0(n+1)}\le \frac{2}{a_0},\quad n\ge 0. 
$$
Thus, \eqref{yeq25} hold if 
\begin{equation}
\label{qeq22}
\frac{2}{a_0} \le \frac{1}{2c_1\lambda}, 
\end{equation}
i.e., if $a_0$ is sufficiently large.

Let us verify inequality \eqref{yeq26}. 
Assume that $(a_n)_{n=0}^\infty$ and $\lambda$ satisfy \eqref{yeq22}--\eqref{yeq25} and \eqref{peq19}. 
Choose $\kappa\ge 1$ such that 
\begin{equation}
\label{eeq14}
\frac{2c_0}{\kappa \lambda} \le \frac{1}{2}.
\end{equation}
Consider the sequence $(b_n)_{n=0}^\infty:=(\kappa a_n)_{n=0}^\infty$ and let $\lambda_\kappa:=\kappa\lambda$. 
Using inequalities \eqref{yeq22}, \eqref{yeq25} and \eqref{eeq14}, one gets
\begin{align*}
c_0\frac{b_n}{\lambda^2_\kappa}+ \frac{b_n-b_{n+1}}{b_{n+1}}c_1 &= 
\frac{2c_0}{\kappa \lambda}\frac{a_n}{2\lambda}+ \frac{a_n-a_{n+1}}{a_{n+1}}c_1\\
&\le  
\frac{1}{2}\frac{a_{n+1}}{\lambda}+ \frac{a_{n+1}}{2\lambda} = \frac{a_{n+1}}{\lambda}= \frac{b_{n+1}}{\lambda_\kappa}.
\end{align*}
Thus, inequality \eqref{yeq26} holds for $a_n$ replaced by $b_n=\kappa a_n$ 
and $\lambda$ replaced by 
$\lambda_\kappa=\kappa\lambda$, where $\kappa \ge \max(1, 
\frac{4c_0}{\lambda})$ (see \eqref{eeq14}). 
Inequalities \eqref{yeq22}--\eqref{yeq25} hold as well under this 
transformation. Thus, the choices $a_n=\frac{a_0\kappa}{n+1}$ and 
$\lambda:= \kappa \frac{M_1}{\|y\|}$, 
$\kappa \ge \max(1, \frac{4c_0\|y\|}{M_1})$, satisfy all the conditions of 
Lemma~\ref{lem1}.
\end{proof}

\begin{rem}{\rm
\label{xrem}
Using similar arguments one can show that the 
choices $\lambda>0$, $a_n=\frac{d_0}{(d+1)^b}$, $d\ge1$, $0<b\le 1,$ 
satisfy all conditions of Lemma~\ref{lem1} 
provided that 
$d_0$ is sufficiently large and $\lambda$ is chosen so that inequality 
\eqref{peq20} holds.
}
\end{rem}

\begin{rem}
\label{rem8}
{\rm
In the proof of Lemma~\ref{lem1} $a_0$ and $\lambda$ can be chosen so 
that $\frac{a_0}{\lambda}$
is uniformly bounded as $\delta \to 0$ regardless of the rate
of growth of the constant $M_1=M_1(R)$ from formula (3) when $R\to\infty$, 
i.e., regardless of the strength of the nonlinearity $F(u)$.
 
Indeed, to satisfy \eqref{peq20} one can choose $\lambda = \frac{M_1}{\|y\|}$.
To satisfy \eqref{qeq21} and \eqref{qeq22} one can choose
$$
a_0 = \max \bigg{(}\sqrt{\lambda \|f_\delta - F(0)\|} , 
4c_1\lambda\bigg{)}\leq \max \bigg{(}\sqrt{\lambda (\|f-F(0)\|+1)} ,
4c_1\lambda\bigg{)},
$$
where we have assumed without loss of generality that $0<\delta< 1$.
With this choice of $a_0$ and $\lambda$, the ratio $\frac {a_0}{\lambda}$
is bounded uniformly with respect to $\delta\in (0,1)$ and does not 
depend on $R$.

Indeed, with the above choice one has $\frac {a_0}{\lambda}\leq 
c(1+\sqrt{\lambda^{-1}})\leq c$,
where $c>0$ is a constant independent of $\delta$, and one can 
assume that $\lambda\geq 1$ without loss of generality.

This Remark is used in the proof of main result in Section 2.2.
Specifically, it will be used to prove that an iterative 
process (25) generates a sequence which stays in a ball
$B(u_0,R)$ for all $n\leq n_0 +1$, where the number $n_0$
is defined by formula (36) (see below), and $R>0$ is sufficiently large.
An upper bound on $R$ is given in the proof of Theorem 10, below formula 
(47).
}
\end{rem}

\begin{rem}
\label{xrem2}
{\rm
It is easy to choose $u_0\in H$ such that 
\begin{equation}
\label{teq20}
g_0:=\|u_0-V_0\|\le \frac{\|F(0)-f_\delta\|}{a_0}.
\end{equation}
Indeed, if, for example, $u_0=0$, then by Lemma~\ref{lem11} and Remark~\ref{remmoi} one gets
$$
g_0=\|V_0\|=\frac{a_0\|V_0\|}{a_0} \le \frac{\|F(0)-f_\delta\|}{a_0}.
$$
If \eqref{yeq23} and \eqref{teq20} hold then
$g_0 \le \frac{a_0}{\lambda}.$
}
\end{rem}

\subsection{Main result}

Let $V_{n,\delta}$ solve the equation:
$$
F(V_{n,\delta}) + a_n V_{n,\delta} - f_\delta = 0.
$$
Denote $V_n:=V_{n,\delta}$. 

Consider the following iterative scheme:
\begin{equation}
\label{3eq12}
\begin{split}
u_{n+1} &= u_n - A_n^{-1}[F(u_n)+a_n u_n - f_\delta],\quad A_n:=F'(u_n)+ a_nI,\quad u_0=u_0,
\end{split}
\end{equation}
where $u_0$ is chosen so that inequality \eqref{teq20} holds.
Note that $F'(u_n)\ge 0$ since $F$ is monotone. Thus, $\|A_n^{-1}\|\le\frac{1}{a_n}$.

Let $a_n$ and $\lambda$ 
satisfy conditions \eqref{yeq22}--\eqref{yeq26}.
Assume that equation $F(u)=f$ has a solution $y\in B(u_0,R)$, possibly nonunique,
and $y$ is the minimal-norm solution to this equation. 
Let $f$ be unknown but $f_\delta$ be given, and $\|f_\delta-f\|\le \delta$.
We have the following result:

\begin{thm}
\label{mainthm}Assume $a_n=\frac{d_0}{(d+n)^b}$ where $d\ge 1,\, 0<b\le 1$, and $d_0$ is sufficiently large
so that conditions \eqref{yeq22}--\eqref{yeq26} hold. 
Let $u_n$ be defined by \eqref{3eq12}. Assume that $u_0$ is chosen so that \eqref{teq20} holds. 
Then there exists a unique $n_\delta$ such that
\begin{equation}
\label{2eq3}
\|F(u_{n_\delta})-f_\delta\|\le C_1\delta^\gamma,\quad
C_1\delta^\gamma < \|F(u_{n})-f_\delta\|,\quad \forall n< n_\delta,
\quad 
\end{equation}
where $C_1>1,\, 0<\gamma\le 1$.

Let $0<(\delta_m)_{m=1}^\infty$ be a sequence such that $\delta_m\to 0$. 
If $N$ is a cluster point of the sequence $n_{\delta_m}$ satisfying \eqref{2eq3}, then
\begin{equation}
\label{feq15}
\lim_{m\to\infty} u_{n_{\delta_m}} = u^*,
\end{equation}
where $u^*$ is a solution to the equation $F(u)=f$.
If 
\begin{equation}
\label{feq16}
\lim_{m\to \infty}n_{\delta_m}=\infty,
\end{equation}
where  $\gamma\in (0,1)$, then
\begin{equation}
\label{feq17}
\lim_{m\to \infty} \|u_{n_{\delta_m}} - y\|=0.
\end{equation}
\end{thm}

\begin{proof}
Denote 
\begin{equation}
\label{beq18}
C:=\frac{C_1+1}{2}.
\end{equation}
Let 
$$
z_n:=u_n-V_n,\quad g_n:=\|z_n\|.
$$ 
We use Taylor's formula and get:
\begin{equation}
\label{1eq9}
F(u_n)-F(V_n)+a_nz_n=A_{a_n} z_n+ K_n, \quad \|K_n\| \le\frac{M_2}{2}\|z_n\|^2,
\end{equation}
where $K_n:=F(u_n)-F(V_n)-F'(u_n)z_n$ and $M_2$ is the constant from \eqref{ceq3}.
From \eqref{3eq12} and \eqref{1eq9} one obtains
\begin{equation}
\label{1eq8}
z_{n+1} = z_n - z_n - A_n^{-1}K(z_n) - (V_{n+1}-V_{n}).
\end{equation}
From \eqref{1eq8}, \eqref{1eq9}, and the estimate $\|A_n^{-1}\|\le\frac{1}{a_n}$, one gets
\begin{equation}
\label{3eq17}
g_{n+1} \le \frac{M_2g_n^2}{2a_n} + \|V_{n+1}-V_n\|.
\end{equation}
Since $0<a_n\searrow 0$, for any fixed $\delta>0$ there exists $n_0$ such that
\begin{equation}
\label{4eq18}
\frac{\delta}{a_{n_0+1}}> \frac{1}{C-1}\|y\|\ge \frac{\delta}{a_{n_0}},\qquad C>1.
\end{equation} 
By \eqref{yeq22}, one has $\frac{a_n}{a_{n+1}}\le 2,\, \forall\, n\ge 0$. This and \eqref{4eq18} imply
\begin{equation}
\label{eeq16}
\frac{2}{C-1}\|y\|\ge \frac{2\delta}{a_{n_0}} >\frac{\delta}{a_{n_0+1}}> \frac{1}{C-1}\|y\|\ge \frac{\delta}{a_{n_0}},\qquad C>1.
\end{equation}
Thus,
\begin{equation}
\label{ceq18}
\frac{2}{C-1}\|y\|> \frac{\delta}{a_{n}},\quad \forall n \le n_0 + 1.
\end{equation}
The number $n_0$, satisfying \eqref{ceq18}, exists and is unique since $a_n>0$ monotonically decays to 0 as $n\to\infty$.
By Lemma~\ref{lem11},
there exists  a number $n_1$ such that 
\begin{equation}
\label{3eq18}
\|F(V_{n_1+1})-f_\delta\|\le C\delta < \|F(V_{n_1})-f_\delta\|, 
\end{equation}
where $V_n$ solves the equation $F(V_{n})+a_{n}V_{n}-f_\delta=0$. 
{\it We claim that $n_1\in[0,n_0]$.} Indeed, 
one has $\|F(V_{n_1})-f_\delta\|=a_{n_1}\|V_{n_1}\|$, and $\|V_{n_1}\|\le \|y\|+\frac{\delta}{a_{n_1}}$ 
(cf. \eqref{2eq1}), so
\begin{equation}
\label{eeq19}
C\delta < a_{n_1}\|V_{n_1}\|\le a_{n_1}\bigg{(}\|y\|+ \frac{\delta}{a_{n_1}}\bigg{)}=a_{n_1}\|y\|+\delta,\quad C>1.
\end{equation}
Therefore,
\begin{equation}
\label{eeq20}
\delta < \frac{a_{n_1}\|y\|}{C-1}.
\end{equation}
Thus, by \eqref{eeq16}, 
\begin{equation}
\label{yeq36}
\frac{\delta}{a_{n_1}} < \frac{\|y\|}{C-1} < \frac{\delta}{a_{n_0+1}}.
\end{equation}
Here the last inequality is a consequence of \eqref{eeq16}.
Since $a_n$ decreases monotonically, inequality \eqref{yeq36} implies $n_1\le n_0$. 
One has
\begin{equation}
\label{eeq21}
\begin{split}
a_{n+1} \|V_n-V_{n+1}\|^2 &= \langle (a_{n+1} - a_{n})  V_{n} - F(V_n) + F(V_{n+1}), V_n - V_{n+1} \rangle \\
&\le \langle (a_{n+1} - a_{n})  V_{n}, V_n - V_{n+1} \rangle \\
&\le (a_{n}-a_{n+1}) \|V_{n}\| \|V_n - V_{n+1}\|.
\end{split}
\end{equation}
By \eqref{2eq1}, $\|V_n\|\le \|y\|+\frac{\delta}{a_n}$, and, by \eqref{ceq18}, 
$\frac{\delta}{a_n}\le \frac{2\|y\|}{C-1}$ for all $ n\le n_0+1$.
Therefore,  
\begin{equation}
\label{ceq40}
\|V_n\|\le \|y\|\bigg{(}1+\frac{2}{C-1}\bigg{)},\qquad \forall n\le n_0+1,
\end{equation}
and, by \eqref{eeq21}, 
\begin{equation}
\label{beq24}
\|V_n-V_{n+1}\| \le \frac{a_n-a_{n+1}}{a_{n+1}}\|V_{n}\|\le \frac{a_n-a_{n+1}}{a_{n+1}}\|y\|
\bigg{(}1+\frac{2}{C-1}\bigg{)},\quad \forall n\le n_0+1.
\end{equation}
Inequalities \eqref{3eq17} and \eqref{beq24} imply
\begin{equation}
\label{1eq10}
g_{n+1}\le \frac{c_0}{a_n}g_n^2+\frac{a_n-a_{n+1}}{a_{n+1}}c_1,
\quad c_0=\frac{M_2}{2},\quad c_1=\|y\|\bigg{(}1+\frac{2}{C-1}\bigg{)},
\end{equation}
for all $n\le n_0+1$. 

By Lemma~\ref{lem1} and Remark~\ref{xrem}, the sequence 
$(a_n)_{n=1}^\infty$, satisfies conditions \eqref{yeq22}--\eqref{yeq26}, 
provided that $d_0$ is sufficiently large and 
$\lambda>0$ is 
chosen so that \eqref{peq20} holds.
Let us show by induction that 
\begin{equation}
\label{ceq15}
g_n<\frac{a_n}{\lambda},\qquad 0\le n\le n_0+1.
\end{equation}
Inequality \eqref{ceq15} holds for $n=0$ by Remark~\ref{xrem2}. Suppose \eqref{ceq15} holds for some $n\ge 0$. 
From \eqref{1eq10}, \eqref{ceq15} and \eqref{yeq26}, one gets
\begin{equation}
\begin{split}
g_{n+1}&\le \frac{c_0}{a_n}\bigg{(}\frac{a_n}{\lambda}\bigg{)}^2 + \frac{a_n-a_{n+1}}{a_{n+1}}c_1\\
&= \frac{c_0 a_n}{\lambda^2} + \frac{a_n-a_{n+1}}{a_{n+1}}c_1\\
&\le \frac{a_{n+1}}{\lambda}.
\end{split}
\end{equation}
Thus, by induction, inequality \eqref{ceq15} holds for all $n$ in the region $0\le n\le n_0+1$.

From Remark~\ref{rem3} one has $\|V_n\| \le \|y\|+\frac{\delta}{a_n}$. 
This and the triangle inequality imply 
\begin{equation}
\label{ceq49}
\|u_0-u_n\| \le \|u_0\|+ \|z_n\|+ \|V_n\|\le \|u_0\|+\|z_n\|+ \|y\|+\frac{\delta}{a_n}.
\end{equation}
Inequalities \eqref{ceq40}, \eqref{ceq15},
and \eqref{ceq49} guarantee that the sequence $u_n$, generated by the 
iterative process \eqref{3eq12}, remains
in the ball $B(u_0,R)$ for all $n\le n_0+1$, where 
$R\le \frac{a_0}{\lambda}+\|u_0\|+\|y\|+ \frac{\delta}{a_n}$.
This inequality and the estimate \eqref{ceq18} imply that the sequence 
$u_n$, $n\le 
n_0+1,$ stays in the ball $B(u_0,R)$,
where 
$$R\le \frac{a_0}{\lambda}+ \|u_0\|+\|y\|+ \|y\|\frac{C+1}{C-1}.$$
By Remark~\ref{rem8}, one can choose $a_0$ and $\lambda$ so that 
$\frac{a_0}{\lambda}$
is uniformly bounded as $\delta \to 0$ even if $M_1(R)\to\infty$ as 
$R\to\infty$ at an arbitrary fast rate.
Thus, the sequence $u_n$ stays in the ball $B(u_0,R)$ for $n\leq n_0+1$
when $\delta\to 0$. An upper bound on
$R$ is given above. It does not depend on $\delta$ as 
$\delta\to 0$.

One has:
\begin{equation}
\label{1eq11}
\begin{split}
\|F(u_n)-f_\delta\|\le& \|F(u_n)-F(V_n)\|+\|F(V_n)-f_\delta\|\\
\le& M_1g_n+\|F(V_n)-f_\delta\|\\
\le& \frac{M_1a_n}{\lambda} + \|F(V_n)-f_\delta\|,\qquad \forall n\le n_0+1,
\end{split}
\end{equation}
where \eqref{ceq15} was used and $M_1$ is the constant from \eqref{ceq3}. 
Since $\|F(V_n)-f_\delta\|$ is nonincreasing, by Lemma~\ref{remark1}, and
$n_1\le n_0$, one gets 
\begin{equation}
\label{3eq21}
\|F(V_{n_0+1})-f_\delta\|\le \|F(V_{n_1+1})-f_\delta\| \le C\delta.
\end{equation}
From \eqref{yeq24}, \eqref{1eq11}, \eqref{3eq21}, 
the relation \eqref{4eq18}, 
and the definition $C_1=2C-1$ (see \eqref{beq18}), one concludes that
\begin{equation}
\label{1eq12}
\begin{split}
\|F(u_{n_0+1})-f_\delta\| 
\le& \frac{M_1a_{n_0+1}}{\lambda} + C\delta \\
\le& \frac{M_1\delta (C-1)}{\lambda\|y\|} + C\delta\le (2C-1)\delta=C_1\delta.
\end{split}
\end{equation}
{\it Thus, if 
$$
\|F(u_0)-f_\delta\|> C_1\delta^\gamma,\quad 0<\gamma\le 1,
$$
then one concludes from \eqref{1eq12} that there exists 
$n_\delta$, $0<n_\delta \le 
n_0+1,$ such that
\begin{equation}
\label{3eq23}
\|F(u_{n_\delta})-f_\delta\| \le C_1\delta^\gamma < \|F(u_{n})-f_\delta\|,\quad 0\le n< n_\delta,
\end{equation}
for any given $\gamma\in (0,1]$, and any fixed $C_1>1$.}

Let us prove \eqref{feq15}. If $n>0$ is fixed, then $u_{\delta,n}$ is a
continuous function of $f_\delta$. Denote
\begin{equation}
\label{ceq46}
\tilde{u}_N=\lim_{\delta\to 0}u_{\delta,N},
\end{equation}
where $N<\infty$ is a cluster point of $n_{\delta_m}$, so that there exists a subsequence of $n_{\delta_m}$,
which we denote by $n_{m}$, such that
$$
\lim_{m\to\infty}n_{m} = N.
$$
From \eqref{ceq46} and the continuity of $F$, one obtains:
$$
\|F(\tilde{u}_N)-f_\delta\| = \lim_{m\to\infty}\|F(u_{n_{\delta_m}})-f_\delta\|\le \lim_{\delta\to 0}C_1\delta^\gamma = 0.
$$
Thus, $\tilde{u}_N$ is a solution to the equation $F(u)=f$, and \eqref{feq15} is proved.

{\it Let us prove \eqref{feq17} assuming that \eqref{feq16} holds.} From \eqref{2eq3} and \eqref{1eq11}  with $n=n_\delta-1$, and from \eqref{3eq23}, one gets
\begin{align*}
C_1\delta^\gamma &\le M_1 \frac{a_{n_\delta-1}}{\lambda} + a_{n_\delta-1}\|V_{n_\delta-1}\|
\le M_1 \frac{a_{n_\delta-1}}{\lambda} + \|y\|a_{n_\delta-1}+\delta.
\end{align*}
If $0<\delta<1$ and $\delta$ is sufficiently small, then
$$
\tilde{C}\delta^\gamma \le a_{n_\delta-1} \bigg{(}\frac{M_1}{\lambda}+\|y\|\bigg{)},\quad \tilde{C}>0,
$$
where $\tilde{C}$ is a constant. 
Therefore, by \eqref{yeq22}, 
\begin{equation}
\label{1eq14}
\lim_{\delta\to 0} \frac{\delta}{2a_{n_\delta}}\le
\lim_{\delta\to 0} \frac{\delta}{a_{n_\delta-1}}\le
\lim_{\delta\to 0} \frac{\delta^{1-\gamma}}{\tilde{C}}\bigg{(}\frac{M_1}{\lambda}+\|y\|\bigg{)}
=0,\quad 0<\gamma<1.
\end{equation}
From \eqref{feq16} and \eqref{1eq14}, by Theorem 6.3.1 in \cite{R499}, one gets \eqref{feq17}.
Theorem~\ref{mainthm} is proved.

\end{proof}


\begin{thebibliography}{99}

\bibitem{D}
K. Deimling, Nonlinear functional analysis, Springer Verlag, Berlin, 1985.

\bibitem{R540}N. S. Hoang and A. G. Ramm, Dynamical systems gradient method for solving ill-conditioned linear 
algebraic systems, (submited).

\bibitem{I}
V. Ivanov, V. Tanana and V. Vasin, Theory of ill-posed problems, VSP, Utrecht, 2002.

\bibitem{M}
V. A. Morozov, Methods of solving incorrectly posed problems, Springer Verlag, New York, 1984.

\bibitem{R499} A. G. Ramm, Dynamical systems method for solving
operator equations, Elsevier, Amsterdam, 2007.


\bibitem{R452} A. G. Ramm, 
Global convergence for ill-posed equations
with monotone operators: the dynamical systems method, J.
Phys A, 36, (2003), L249-L254.

\bibitem{R454} A. G. Ramm,
 Dynamical systems method for solving nonlinear
operator equations, International Jour. of
Applied Math. Sci., 1, N1, (2004), 97-110.

\bibitem{R456} A. G. Ramm,
Dynamical systems method for solving
operator equations, Communic. in Nonlinear Sci. and
Numer. Simulation, 9, N2, (2004), 383-402.


\bibitem{R469} A. G. Ramm,
DSM for ill-posed equations with monotone 
operators,
Comm. in Nonlinear Sci. and Numer. Simulation, 10, N8, (2005),935-940.

\bibitem{R474} A. G. Ramm,
Discrepancy principle for the dynamical 
systems
method, Communic. in Nonlinear Sci. and Numer. Simulation,
10, N1, (2005), 95-101

\bibitem{R485} A. G. Ramm, 
Dynamical systems method (DSM) and
nonlinear problems, in the book: Spectral Theory and Nonlinear 
Analysis,
World Scientific Publishers, Singapore, 2005, 201-228. (ed J.
Lopez-Gomez).


\bibitem{R491} A. G. Ramm,
Dynamical systems method (DSM) for unbounded
operators, Proc.Amer. Math. Soc., 134, N4, (2006), 1059-1063.




\end{thebibliography}
\end{document}